\documentclass[11pt]{article}
\usepackage{graphicx}
\usepackage{latexsym}
\usepackage{amsmath}
\usepackage{amssymb}
\def\dsum{\displaystyle\sum}
\begin{document}

\title{$\ \ $The structure of HCMU metric in a K-Surface}
\author{ Qing Chen \ \ \,   Xiuxiong Chen \ \ \ and \ \  Yingyi Wu }




\date{}%
\maketitle

\begin{abstract}We study the basic structure of a HCMU metric in a
K-Surface with prescribed singularities. When the underlying
smooth surface is $S^2$, we prove the necessary condition given in
[1] for the existence of HCMU metric is also sufficient.

\end{abstract}
\maketitle
\section{Introduction}

Let M be any compact, oriented smooth Riemannian surface without
boundary, and $M_{\{\alpha_1, \alpha_2, \cdots, \alpha_n\}}$ (
where $\alpha_i > 0, \forall i,\  1 \le i \le n$ ) denotes a
K-Surface associated with M. A Riemannian metric $g$ is said to be
well defined or smooth in $M_{\{\alpha_1, \alpha_2, \cdots,
\alpha_n\}}$ if it satisfies the following two conditions:
\begin{description}\item[](1)\ \ $g$ is smooth everywhere on M except in a set of
singular points $\{ p_1, p_2, \cdots, p_n \}$, \item[](2)\ For any
$i (1 \le i \le n)$, the metric $g$ has a singular angle of $2\pi
\alpha_i$ at the point $p_i$.\end{description} Here the condition
(2) means that in a small neighborhood of $p_i$, there exists a
local complex coordinate chart $ ( U,z)$ $ ( z(p_i)=0)),$ s.t.
$$ g |_U=h(z,\bar{z}){1 \over {{|z|}^{2-2{\alpha}_i}}}{|dz|}^2
,$$
 where $h:~U\to R$ is a continuous positive function and smooth on $U\setminus \{0\}$.
  Two smooth
Riemannian metrics on $M_{\{\alpha_1,\alpha_2,\cdots,\alpha_n\}}$
are pointwise conformal to each other if they are related by a
multiple of a smooth positive function on $M$.

A natural question is whether or not there exists a ``best" metric
in every conformal class of a K-Surface. This is an attempt to
generalize  the classical uniformization theorem to  a K-Surface.
Recalled that the classical uniformization theorem asserts that in
every conformal class of $M$, there must exist a metric with
constant scalar curvature. Many papers tried to generalize the
uniformization theorem in  K-Surfaces. For example, [5] and [3]
independently found the sufficient condition under which in a
K-Surface, there exists a constant scalar curvature metric. [6]
found a necessary condition of the existence of a constant scalar
curvature metric in a K-Surface, [9] proved that a uniqueness
theorem on constant curvature metric in some K-surfaces. However,
there does not exist a constant curvature metric in a K-Surface.

In a serial of papers [1] and [2], the second named author tried
to find the ``best" metric
 in a conformal class of a K-Surface,
  through studying the critical point of the Calabi energy functional. He proposed that
two types of metrics can be regarded as candidates of the ``best"
metric: one is the extremal metric, another is the HCMU metric (
Definitions will be given later.).

Let $M_{\{\alpha_1,\alpha_2,\cdots,\alpha_n\}}$ be a K-Surface and
$g_0$ be a smooth metric in it. Consider the conformal class of
$g_0$ : $$\mathcal{S}(g_0)=\{ g= e^{2 \varphi}g_0,\varphi \in
H^{2, 2}(M) \mid \int_{M\setminus \{p_1, p_2, \cdots, p_n \}} e^{2
\varphi}dg_0=\int_{M\setminus \{p_1, p_2, \cdots, p_n \}}
dg_0\}.$$Define the Calabi energy functional:\begin{equation}
E(g)=\int_{M \setminus \{p_1,p_2,\cdots,p_n\}} K^2
dg,\end{equation} here $K$ is the scalar curvature of the metric
$g$. The Euler-Lagrange equation of $E(g)$ is (cf. [1] [8] )
\begin{equation} \triangle_g
K+K^2=C,\end{equation}
 or equivalently, in a local complex coordinate chart,
 \begin{equation} {\partial \over
{\partial \bar{z}}} K_{,zz}=0,\end{equation} where $K_{,zz}$ is
the 2nd-order $(0,2)$ type covariant derivatives of $K$.

 A metric which
satisfies (2) or (3) is called an extremal metric. (3) has two
special cases, one is
\begin{equation} K\equiv Const,
\end{equation} and the other is \begin{equation}
K_{,zz}=0, \ K \neq Const. \end{equation} A metric which satisfies
(5) is called a HCMU (the Hessian of the Curvature of the Metric
is Umbilical) metric. Throughout this paper, we assume that a HCMU
metric has finite area and finite Calabi energy.

Let us first quote an Obstruction Theorem from [1].\par

 {\bf Theorem 1.} {\it Let $g$ be a HCMU metric in a K-Surface
    $M_{\{\alpha_1, \alpha_2, \cdots, \alpha_n\}}$. Then the Euler
    character of the underlying surface should be determined by
    \begin{equation}\chi(M)=\sum_{i=1}^j(1-\alpha_i)+(n-j)+s\end{equation}
    where $s$ is the number of critical points of the Curvature
    $K$ ( excluding the singular points of $g$). Here we assume
     that $\alpha_1, \alpha_2, \cdots, \alpha_k, ( 0\le k \le n)$  are
     the only integers in the set of prescribed angles
     $\{\alpha_1, \alpha_2, \cdots, \alpha_n\}$; and assume that
     $\{p_{j+1}, p_{j+2}, \cdots, p_k\}$ are the only local extremal
     points of $K$ in the set of singular points $\{p_j, 0\le j
     \le k \}$.\par}

The formula (6) is an application of Poincar$\grave{\mbox{e} }$
-Hopf index theorem. When $g$ is a HCMU metric, the gradient
vector field $\overrightarrow{V}$ of the scalar curvature $K$ is
holomorphic. Hence, its real part is a Killing vector field. It
was proved in [1] that the singularities of the Killing vector
field is a finite set which is the union of the singularities of
metric $g$ and the smooth critical points of function $K$.
Consequently, any saddle point of $K$ must be the singularities of
metric $g$. At these points the index of the vector field is
$(1-\alpha_i)$. Other singularities of this gradient vector field
must be local extremal points of $K$ with index $1$. Therefore,
the Poincar$\grave{\mbox{e}}$-Hopf index theorem implies formula
(6).

 In this paper, we study the following question: whether or not
the condition (6) is also sufficient to the existence of  HCMU
metrics in a K-Surface. Our main result in this paper is:

      \noindent{\bf Theorem  A.} For $S^2$,  given $n$ points
     $p_1,p_2,\cdots,p_n$ on $S^2$ and $n$ positive numbers $\alpha_1,
       \alpha_2, \cdots, \alpha_n
       $ with $ \alpha_1, \alpha_2, \cdots, \alpha_k$ being the only integers and $ \alpha_j \ge 2 (1 \le j \le k)$, suppose that $\alpha_1, \alpha_2, \cdots,
       \alpha_n$
     satisfy the following condition: $ \exists j_0 (1 \le j_0 \le k)$   and  $\alpha_{\sigma{(1)}}, \alpha_{\sigma{(2)}},
      \cdots, \alpha_{\sigma{(j_0)}}, \sigma(i)\in\{1, 2, \cdots
     k\}(1 \le i \le j_0)$, s.t.\begin{equation}\sum_{i=1}^{j_0}\alpha_{\sigma(i)}+\chi(M)-n \ge 0.\end{equation}Then  there exists a
     HCMU  metric whose scalar curvature $K$ is not a constant,  s.t.
     the angles of the metric at $p_1,p_2,\cdots,p_n$ are exactly
     $\alpha_1, \alpha_2, \cdots, \alpha_n$ and $p_{\sigma{(1)}}, p_{\sigma{(2)}}, \cdots, p_{\sigma{(j_0)}}
     $ are the only saddle points of the scalar curvature $K$.

     In fact, we prove for $S^2$ the condition (7) is the necessary and sufficient
     condition for the existence of a HCMU metric in it.

     The simplest HCMU metric in $S^2$ is a football. It only has
     two extremal points and it is a rotationally symmetric
     metric. Fixing the area, a HCMU metric is uniquely determined by
     the ratio of the two angles.

     The proof of Theorem A is based on the following Theorem B, which says that any HCMU metric can be divided into a
     finite number of
     footballs. Under the condition (7), we can glue some suitable footballs together to
     obtain a HCMU metric in $S^2$ as desired.

     \begin{flushleft}{\bf Theorem B.} Let $g$ be a HCMU metric on a K-Surface
     $M$, then there are a finite number of geodesics which connects extremal points and saddle points of the scalar
      curvature $K$ together. In fact, $M$
can be divided into a finite number of pieces by cutting along
these
      geodesics where each piece is locally isometric to a HCMU
      metric in some football.\end{flushleft}

     We should point out that in [7], Lin and Zhu use ODE method and geometry of the scalar curvature of HCMU
     metrics to construct a class of HCMU metrics with
     finite conical singular angles $2 \pi\cdot integers$ on
     $S^2$. This kind of HCMU metric is called
     exceptional HCMU metric where all of its singularities are the saddle points of the scalar curvature $K$. A minimal exceptional
     HCMU metric is an exceptional HCMU metric with only one minimum point of the scalar curvature $K$.
      They give an explicit formula for minimal exceptional HCMU
     metrics. Their theorem shows that a minimal exceptional HCMU metric is determined
     by three parameters. In comparison, our existence theorem of HCMU metric is more general. Indeed
     our construction in the proof of Theorem A is actually a minimal
     exceptional HCMU metrics if all of the singularities are the
     saddle points of the scalar curvature $K$.

     The authors would like to thank the referee for many useful
     suggestions.
\section{Proof of Theorem B}
\subsection{Preliminaries}

Let $M$ be a compact, oriented smooth Riemannian surface without
boundary. $M_{\{\alpha_1,\alpha_2,\cdots,\alpha_n\}}$ denotes its
K-Surface. $\{p_1,p_2,\cdots,p_n\}$ is the set of singular points.
$\{\alpha_1,\alpha_2,\cdots,\alpha_n\}$ is the corresponding set
of singular angles. $\forall p \in M \setminus
\{p_1,p_2,\cdots,p_n\}$, assuming that $(U,z)$ is a complex
coordinate chart around $p$, $g$ can be written as: $$
g=e^{2\varphi(z,\bar{z})}{|dz|}^2.$$ and $$K=-{\triangle \varphi
\over {e^{2\varphi}}}.$$ Equation (5) can be written
as:\begin{equation} K_{, zz}={\partial^2 K \over {\partial z}
^2}-2
   {\partial K \over \partial z}{\partial \varphi \over  \partial
   z}=0,\end{equation}
which means that the gradient vector field of the scalar curvature
$K$ is holomorphic. The gradient vector field ${\nabla K}$ is:
$$ {\nabla K}=\sqrt{-1}{K}^{\prime
   \bar{z}}{\partial \over \partial z}=\sqrt{-1} e^{-2 \varphi}{\partial K \over \partial \bar{z}}{\partial
   \over \partial
   z},$$
   and its real part is: $$\overrightarrow {V}={1 \over
2}(\sqrt{-1}{K}^{\prime z}  {\partial \over \partial
z}-\sqrt{-1}{K}^{\prime \bar{z}}{\partial \over \partial \bar
{z}}).$$ Then $\overrightarrow {V}$ is a Killing vector field and
its integral curve is the level set of the function $K$.

In fact, by studying the properties of the Killing vector field
$\overrightarrow{V}$, the Obstruction Theorem is proved in [1]. We
list here the main properties of $\overrightarrow {V}$.

\par\vskip0.3cm
\noindent {\bf Proposition 1[1].} Let $Sing \overrightarrow{V} $
denote the set of all singular points of $\overrightarrow{V}$ and
$\Omega_p$ denote the set of the integral curves of
$\overrightarrow{V}$ which meet $p$ for $p \in Sing
\overrightarrow {V}$. Then:\\
(1) $Sing \overrightarrow{V} = \{ \ smooth \ critical \ points \
of \ K \} \bigcup \{p_1,p_2,\cdots,p_n\}$ and is a finite
set.\\
(2) $\Omega_p$ is empty or a finite set. Moreover, if
$\Omega_p\neq \emptyset$, $\Omega_p$ has even number of points and
$g$ has angle $|\Omega_p|\pi$ at $p$. \\
(3) $K$ can be continuously extended to $M$.\\
(4) $Sing \overrightarrow{V}$ can be divided
 into two parts, $ Sing \overrightarrow{V}=S_1 \bigcup S_2,\ $
 such that\\
$~~~$ (a) $S_1=\{p \in Sing \overrightarrow{V}
 |\Omega_p=\emptyset \}$,\ and if $p \in S_1$, then $p$ is an extremal point of
 $K$;\\
 $~~~$ (b) $S_2 =\{p \in Sing \overrightarrow{V}|\Omega_p \neq \emptyset\}$, if
 $p \in S_2$, $p$ is a saddle point of $K$. and at $p$ the angle is
 $\pi\cdot|\Omega_p|.$ 

\par\vskip0.3cm
 \noindent{\bf Proposition 2[1].}\\ 
 (1) Any  integral curve of ${\overrightarrow
     {V}}$ in the neighborhood of any local extremal of $K$ point is a topologically  circle
     which contains the point in its interior. \\
     (2) If a closed
     integral curve of $\overrightarrow{V}$ bounds a
     topological disk in $M$, which contains only one extremal
     point of $K$, then every integral curve of
     $\overrightarrow{V}$ in this disk is  also a topological
     circle. 

     \par\vskip0.3cm
     \noindent{\bf Proposition 3[1].} At a saddle point of $K$, the
     included angle of two adjacent integral curves of $\nabla K$ is $\pi$.

     \subsection{Proof of the Theorem B}

     Let us begin with the study of a football, i.e. a HCMU metric $g$ in $S^2$, which is
     rationally symmetric and  has two extremal points.

     According to Proposition 1, the scalar curvature $K$ is continuous. If
     $K$  has only two extremal points $p$ and $q$, by Proposition 2.1, in a neighborhood of $p$, an integral curve $C$ of $\overrightarrow
     {V}$ bounds a topological disk $D$ centered at $p$. By of Proposition 2.2, the
      integral curves of $\overrightarrow
     {V}$ in $D$ are all
     topologically concentric circles containing $p$ in their interiors. Since
     $\overrightarrow {V}$ is a Killing
     vector field, $g$ is invariant along integral curves of $\overrightarrow {V}$, then
      $g$ is rotationally symmetric in $D$.
      On the other hand,  $C$ is also a topological circle bounding
      the disk $S^2 \setminus D$ which has only one extremal point $q$, by
      Proposition 2.2 again, $g$ is also rotationally symmetric in $S^2 \setminus D$. Therefore, $g$ is globally rotationally symmetric.
     It can be written as:
     \begin{equation}g={du}^2+f^2(u) {d\theta}^2
     \ \ (0 \le u \le l,\ 0 \le \theta \le 2 \pi),\end{equation}
     with $p$ and $q$ corresponding to $u=0$ and $l$
     respectively, and $\mbox{dist}_g(p,q)=l,$
     see Figure 1.

     \begin{center}

     \includegraphics{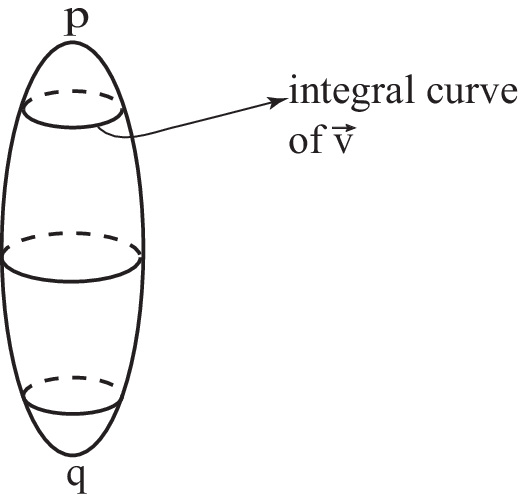}\\
     Figure 1
     \end{center}

     If we assume the angle of the
     metric at $p$ is $\alpha$, the angle at $q$ is $\beta$, $(\alpha
     \ge \beta )$, then $f$ satisfies:

    \parbox{9.75cm}{\[\left \{ \begin{array}{l}
     f(0) = f(l)=0, \\ f^\prime(0)=\alpha, f^\prime(l)=-\beta, \\
    f(u)>0 ,  u  \in (0, l). \end{array}\right.\]}\hfill
    \parbox{1cm}{\begin{eqnarray}\end{eqnarray}}

    \noindent By (9), the scalar curvature $K$ is given by:\begin{equation}K=-{f^{\prime \prime} \over
    f}.\end{equation}
    \par
\vskip0.3cm
    \noindent {\bf Proposition 4.} There is a constant $c$ such that ${K}^\prime=cf
    .$\par\vskip0.1cm
    \noindent{\it Proof:} Since $g$ is a HCMU
    metric, from (8) and (9) we have ${K}^{\prime \prime} f={K}^\prime f^\prime,$
    that is
    ${K}^{\prime}=cf.$ $\Box$

    \par\vskip0.3cm
    \noindent {\bf Proposition 5.} ${K}^\prime \le 0, $ moreover, if
    $K \neq Constant, $ then only when $u=0 \  \mbox{and}
    \ u=l,\  {K}^\prime =0.$\par\vskip0.1cm
    \noindent{\it Proof:} Define a function: $F={f^\prime}^2+K f^2.$ Then
    $F(0)=\alpha^2,\  F(l)=\beta^2,\  F(l) \le F(0), $ so there is a  $\xi \in (0, l)$ such that
    $${F(l)-F(0) \over
    l}=F^\prime (\xi)\le 0.$$
    On the other
    hand,
    $$F^\prime=2f^\prime f^{\prime \prime}+{K}^\prime
    f^2+2K f f^\prime,$$
     by $K=-{f^{\prime \prime}\over
    f},$ we have ${K}^\prime (\xi) \le 0.$ From {Proposition
    4}, ${K}^\prime =cf$, and $f$ is positive on $(0,l)$, if $K$ is not a constant,
    $K^\prime$ does not change its sign from $0$ to $l$. Hence, $K^\prime
    (\xi) \le 0$ implies ${K}^\prime(u) < 0,\ (\forall u \in (0, l))$. Moreover,
    $K^\prime=0$ only when $u=0$ and $l$.     $\Box$

    \par\vskip0.3cm
    \noindent {\bf Remark 1:} From the proof of {Proposition 5}, we also get $K
    \equiv Constant$ if and only if $\alpha=\beta$.
\par\vskip0.3cm
    In the following, we always assume $K \ne
    Constant$. Therefore $K$ decreases monotonely from $p$ to
    $q$.  Substituting $f=\frac{{K}^\prime}c $
    into $K=-{f^{\prime \prime} \over f},$ we
     get $K'''+K'K=0$, that is
     \begin{equation} {{{K}^\prime}^2 \over
     2}=C_0K-{{K}^3 \over 6} +C_1,\end{equation}
     here $C_0$ and $C_1$ are two constants. Assuming $K(0)=K_0, K(l)=K_1$ and
     letting $u=0$ and $l$ in
     (12), we know both $K_0$
      and $K_1$ are  roots of the equation $-{{K}^3
      \over 6}+C_0K+C_1=0$. Then
       $$-{{K}^3\over
      6}+C_0K+C_1=-{1 \over 6}(K-K_0)(K-K_1)(K+K_0+K_1)$$
      and \begin{equation} {{K}^\prime}^2=-{1 \over 3}
      (K-K_0)(K-K_1)(K+K_0+K_1).\end{equation}
      We take derivatives of (13) to
      get:\begin{equation}{K}^{\prime \prime}=-{1 \over
      6}[(K-K_0)(K+K_0+K_1)+(K-K_1)(K+K_0+K_1)+(K-K_0)(K-K_1)].\end{equation}
      Using {Proposition 4}, we have:\begin{equation}
      c f^\prime=-{1 \over
      6}[(K-K_0)(K+K_0+K_1)+(K-K_1)(K+K_0+K_1)+(K-K_0)(K-K_1)].\end{equation}
      Let $u=0 \ \mbox{and}\
      l$ in (15), then by (10):\[ \left \{ \begin{array}{ll}
      \alpha=\dfrac{(K_1-K_0)(2K_0+K_1)}{6c},
      \\\beta=\dfrac{(K_1-K_0)(2K_1+K_0)}{6c}.\end{array}\right.
      \]
      Integrate the equation ${K}^\prime=cf$ from $0$ to $l$, we get:\begin{equation}
      K_1-K_0=c\int_0^lf(u)du=c{A(g) \over 2\pi}, \end{equation}
      here $A(g)$ denotes the area of the metric $g$. Then we have:

      \begin{equation} \left\{
      \begin{array}{l}\alpha=\dfrac{A(g)}{12 \pi} (2K_0+K_1), \\
      \beta=\dfrac{A(g)}{12 \pi}(2K_1+K_0).
      \end{array} \right.
    \end{equation}

      or:

\begin{equation}  \left \{ \begin{array}{l}
K_0=\dfrac{4\pi}{A(g)}(2\alpha-\beta), \\ K_1=\dfrac{4\pi}{A(g)}
      (2\beta-\alpha).\end{array} \right.
      \end{equation}

      From (18) we see
      if $\alpha,\  \beta,\  A(g)$ are fixed, then $K_0$ and $K_1$ are uniquely
      determined, and  by (13) $K$ is determined, again according to ${K}^\prime=cf,$ $f$ is
      determined, i.e. the metric $g$ is determined. Therefore we get:

      \begin{flushleft}{\bf Theorem C.} If area and angles at both extremal points are given, there exists a unique
      rotationally symmetric HCMU metric in $S^2$, which is a football.\end{flushleft}

      Meanwhile we get:

      \begin{flushleft}{\bf Corollary 1.} In a football, assume that $\alpha$ is the angle of the HCMU
      metric at the local maximum point of $K$, $\beta$ is the angle of the metric at the local minimum point of
       $K$, $K_0=\max K$, $K_1=\min K$,  then \\
      (1) $K_0>0$, the sign of $K_1$ is the same as
      $2\beta-\alpha.$\\
      (2)  $K_0>K_1>-(K_0+K_1).$\end{flushleft}

      Next we consider a HCMU metric $g$ in a K-Surface $M_{\{\alpha_1,\alpha_2,\cdots,\alpha_n\}}$. Since the integral curves of $\overrightarrow {V}$ are the
      level sets of $K,\  \nabla K \bot \overrightarrow
      {V}$, integral curves of $\nabla K$ are geodesics. If $p$ is
      a local minimum point of $K$,  in a small neighborhood of $p$, the integral curves
       of $\overrightarrow {V}$ are topologically concentric circles and the integral curves of $\nabla K$ are perpendicular to
        them. Choose an integral curve $c(t)\ ~(t\in [0,T])$ of $\overrightarrow {V}$. For each $t$, there
        exists a unique integral curve $C_t$ of $\nabla K$ starting from $p$ and
        passing through the point $c(t)$. See Figure 2.

       \begin{center}
       \includegraphics{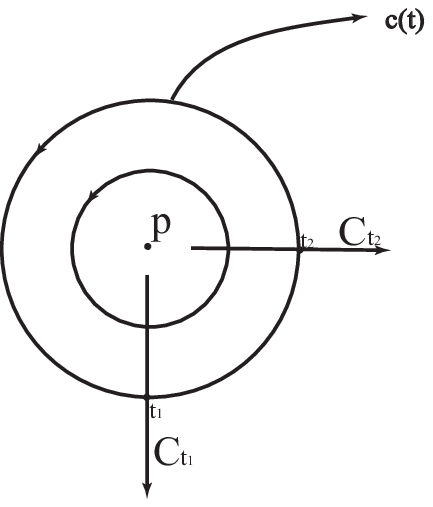}\\
       Figure 2
       \end{center}

        Obviously, $C_t$ must reach some saddle point of $K$ or some local
        maximum point of $K$. We have the
        following:\par\vskip0.3cm
        \noindent{\bf Lemma 1.} For $t_1,t_2 \in (0,t_0)$, If both $C_{t_1}$ and $C_{t_2}$
        reach local
        maximum points $q_1$ and $q_2$ directly without passing through any
         saddle point of $K$, then
         $\mbox{dist}_g(p,q_1)=\mbox{dist}_g(p,q_2)=l$.\par\vskip0.2cm
        \noindent{\it Proof:} Notice that in a small neighborhood of $p$ (we say
         $D$), $g$ is rotationally symmetric. Therefore, in $D$, $g$ can be
         written as $g={du}^2+f^2(u){dt}^2,$ where $t$ is the parameter of $c(t)$. Hence
         $K=-{f^{\prime \prime} \over f},\ {K}^\prime=\frac{dK}{du}=cf.$
         Then by (13)
         \begin{equation} {K}^\prime=-\sqrt{-{1 \over 3} (K-K_0)({K}^2+a_0
         K+a_1)}, \end{equation}\\ here $K_0=K(p),\ a_0 \ \mbox{and} \ a_1 \ \mbox{are
         constants}$. On the other hand, if we restrict
         $K$ at $C_{t_1}$,  $K$ is a smooth function of the arc
         length parameter $s\ (s=u)$ of $C_{t_1}$. Moreover, at $C_{t_1} \bigcap D$, we have
         (19). According to ODE theory, we  know that (19) holds true at the whole of $C_{t_1}$.
          Since ${K}^\prime
         (q_1)=0$ (It is because in a neighborhood of $q_1$, $g$ is also rationally
         symmetric, we  have $K'=\tilde{c}\tilde{f}$ and
         $\tilde{f}(0)=0$),
          we get:\begin{equation} {K}^\prime=-\sqrt{-{1 \over
         3}(K-K_0)(K-K_1)(K+K_0+K_1)}, \end{equation}here $K_1=K(q_1)$. At $C_{t_2}\bigcap
         D,$ we also have
         $$\ {K}^\prime =-\sqrt{-{1 \over
         3}(K-K_0)(K-K_1)(K+K_0+K_1)}.$$ Because along
         $C_{t_2}$ except at $p$ and $q_2$, there is no point at which $\nabla K=0$, we get
         $K(q_2)=K_1$. Furthermore, since
         $${dK \over ds}=-\sqrt{-{1 \over
         3}(K-K_0)(K-K_1)(K+K_0+K_1)},$$ \begin{equation}{ds \over dK}=-\dfrac{1}{ \sqrt{-{1 \over
         3}(K-K_0)(K-K_1)(K+K_0+K_1)}}. \end{equation} Hence, the length of
         a geodesic from $p$ to $q_1$ is:\begin{equation}
         l=\int_{K_1}^{K_2}-{dK \over \sqrt{-{1 \over 3}
         (K-K_1)(K-K_2)(K+K_1+K_2)}}. \end{equation}The same
         as above,  the length of a geodesic from $p$ to $q_2$ is
         also $l$.$\Box$

         \par\vskip0.3cm
         \noindent{\bf Lemma 2.} Fix $t_0 \in [0,T]$ and suppose
         $C_{t_0}$ reaches a maximum point $q$ of $K$ directly,
         then $\exists \varepsilon>0$, s.t. $\forall t \in
         (t_0-\varepsilon,t_0+\varepsilon)$, $C_t$ reaches the same
          maximum point $q$ without passing through any saddle point of
          $K$.\par\vskip0.2cm
        \noindent{\it Proof:} Since the saddle points of $K$ are finite,
        there exists a small neighborhood of $t_0 \ (t_0-\varepsilon
        ,t_0+\varepsilon)$ s.t. each $C_t (t \in (t_0-\varepsilon
        ,t_0+\varepsilon))$ does not reach any saddle point. The end points of $C_t$
        are continuously dependent on $t$
         and the maximum points of $K$ are finite. Therefore, $\exists \varepsilon>0$ s.t.
         $\forall t \in (t_0-\varepsilon,
         t_0+\varepsilon)$, $C_t$ reaches the same maximum point.
          $\ \Box$
\par\vskip0.3cm
         \noindent{\it Proof of Theorem B:} By virtue of Lemma 2,
         $\bigcup\limits_{t \in (t_0-\varepsilon,t_0+\varepsilon)}C_t$ is a simply
         connected domain in $M$. Suppose $F$ is the largest simply connected domain in $M$ which contains $\bigcup\limits_{t \in (t_0-\varepsilon,t_0+\varepsilon)}C_t$
           and satisfies the following three properties:
\begin{enumerate}\item[1] Any integral curve of $\nabla K$ in $F$
           is from $p$ to $q$. \item[2] Any
          integral curve of $\nabla K$ in $F$ does not pass through any
         saddle point of $K$. \item[3] $\partial F$ are also integral curves of $\nabla K$, but they pass through some
         saddle points of $K$.\end{enumerate}

         $\partial F$ can be divided into two curves connecting
         $p$ and $q$, we say $\gamma_1$ and $\gamma_2$. Along
         $\gamma_1$, there are  some saddle points which are connected by
         geodesic segments (integral curves of $\nabla K$). At
         each saddle point, by Proposition 3, the included angle
         of two adjacent geodesics is $\pi$. Hence, $\gamma_1$ is
         smooth at each saddle point. Therefore, $\gamma_1$ is a
         smooth geodesic, either is $\gamma_2$.\par
         On the other hand, in $F$, $g$ is invariant along integral curves of $\overrightarrow
         {V}$, i.e. $g$ is rotationally symmetric in $F$.
         Therefore, we can parameterize $g$ as: $
         g={du}^2+f^2(u){d\theta}^2 $ as before. Hence, $(F,g)$ is
         isometric to a football.\par
         At a minimum point $p$ of $K$, since there are finite integral curves of $\nabla K$
         starting from $p$ and reaching saddle points, we can
         repeat the above proof to obtain finite pieces of the
         largest simply connected domains in $M$, which are all
         isometric to footballs and contain $p$ as a vertex. We
         can also repeat this operation at each minimum point of
         $K$ to obtain finite footballs. We claim that every point
         of $M$ is contained in the union of these footballs. For
         any point of $M$, there must be an integral curve of
         $\nabla K$ starting from it or ending at it or passing
         through it. If there are some saddle points at this integral curve, the point must
         be on the boundary of some largest domain. If not, this
         point must be in some largest domain. Therefore, every
         point of $M$ is contained in the union of these
         footballs. This completes the proof of Theorem B.
         $\Box$\par

     We also have following corollaries.

     \begin{flushleft}{\bf Corollary 2. }On HCMU surface, the values of local minimum of scalar
     curvature $K$  are the same to each other, and the same narration is true
     for local maximum points of $K$.\end{flushleft}
     {\bf Corollary 3.} Let $g$ be a HCMU metric in a K-surface, If scalar curvature $K$
     has only extremal  points in the Surface, then the
     metric is a football.

      \section{Proof of Theorem A}

  We prove following lemma at first.\par\vskip0.2cm
  \noindent {\bf Lemma 3.} Two footballs $S^2_{\{\alpha,\beta\}}$ and $S^2_{\{\alpha_1,\beta_1\}}$ can be smoothly glued along their meridians or some
  segments of the two meridians, iff ${\alpha \over \beta}={\alpha_1
  \over \beta_1}$ and $\frac{A(g)}{A(g_1)}=\frac{\alpha}{\alpha_1}$, here
  $A(g)$ and $A(g_1)$ denote the areas of the two
  metrics.\\ \vskip0.05cm
  {\it Proof:} Suppose $K_0~(\tilde{K}_0)$ and $K_1~(\tilde{K}_1)$ are the maximum and minimum of
  $K~(\tilde{K})$ of football $S^2_{\{\alpha,\beta\}}~(
  S^2_{\{\alpha_1,\beta_1\}})$ respectively.\par
  $(\Longrightarrow)$ If two footballs can be smoothly glued along their meridians or some
   segments of the two meridians, we see at each segment being smoothly glued, the arc length is
    the
  same, $K$ is the same, so the derivative of $K$, ${K}^\prime$ is the same. Therefore
  we get from the equation (13) and {Corollary 1}, $K_0=\tilde{K_0},\ K_1=\tilde{K_1}.$
  Then from
  the equation (17), we get ${\alpha \over \beta}={\alpha_1 \over
  \beta_1}$ and $\frac{A(g)}{A(g_1)}=\frac{\alpha}{\alpha_1}.$\par
  \vskip0.1cm
   $(\Longleftarrow)$ Assume the metric of two football are given by
   \begin{eqnarray*}
  g &=& {du}^2+{\alpha_1}^2 f^2
  {d\theta}^2 \ (0 \le u \le l, 0 \le \theta \le \frac{2 \pi}{\alpha_1}),\\
  \tilde{g} &=& {du_1}^2+\alpha^2 {f_1}^2{d\theta_1}^2 \ (0 \le u_1 \le l_1,0 \le \theta_1 \le \frac{2 \pi}{\alpha}).
  \end{eqnarray*}
  If $\frac{\alpha}{\beta}=\frac{\alpha_1}{\beta_1}$ and
  $\frac{A(g)}{A(g_1)}=\frac{\alpha}{\alpha_1}$, by the equation
  (18), we get $K_0=\tilde{K_0},K_1=\tilde{K_1}$. Then by the
  equation (21) and (20), $u_1=u,K'=\tilde{K}'$ and by the
  equation (22), $l=l_1$. Meanwhile, we have
  $\frac{f}{f_1}=\frac{\alpha}{\alpha_1}.$ Hence, if we let
  $u=u_1$ and $\theta=\theta_1$, then $\alpha_1 f=\alpha f_1$, i.e. $g$ and $\tilde{g}$ are locally
  isometric. Since the  meridians are geodesics, two footballs
   $S^2_{\{\alpha,\beta\}}$ and $S^2_{\{\alpha_1,\beta_1\}}$ can
   be glued together along their meridians or some segments of
   meridians.\ $\Box$\par\vskip0.2cm

  \noindent{\it Proof of Theorem A.}  ~~\\

   {\it Claim: there exists a HCMU metric on $S^2$ whose singular angles
      are exactly $\alpha_1,\alpha_2,\cdots,\alpha_n$ and the
      angles at the $j_0$ saddle points of the scalar curvature $K$ are exactly $\alpha_{\sigma(1)},\alpha_{\sigma(2)},\cdots,
      \alpha_{\sigma(j_0)}$.}\\

      We will use some suitable footballs to construct $S^2_{\{\alpha_1, \alpha_2, \cdots, \alpha_n
  \}}$ which satisfies the condition of the claim. Without loss of generality, we assume $\alpha_{\sigma(i)}=\alpha_i,\
  i=1,2,\cdots,j_0.$ \par\vskip0.05cm
  \noindent{\it Step 1: construction of a HCMU metric with one saddle point of
     angle $2 \pi \alpha_i \ (i=1,2,\cdots,j_0).$}\par
     Choose $\alpha_i$ footballs, say
     $S^2_{\{x_1,y_1\}},S^2_{\{x_2,y_2\}},\cdots,S^2_{\{x_{\alpha_i},y_{\alpha_i}\}},$
     they satisfy  $$\frac{x_m}{y_m}=\frac{x_n}{y_n}>1,~~~
     \frac{A(g_m)}{A(g_n)}=\frac{x_m}{x_n},$$
     for any $m,n \in
     \{1,2,\cdots,\alpha_i\}$. Then by Lemma 3, two
     footballs can be glued together. We are going to glue these
     footballs together. For example, when $\alpha_i=3,$ take 3
     footballs
     $S^2_{\{x_1,y_1\}},~~~~S^2_{\{x_2,y_2\}}~$, $S^2_{\{x_3,y_3\}}$, see
     Figure 3.
  \begin{center}
       \includegraphics{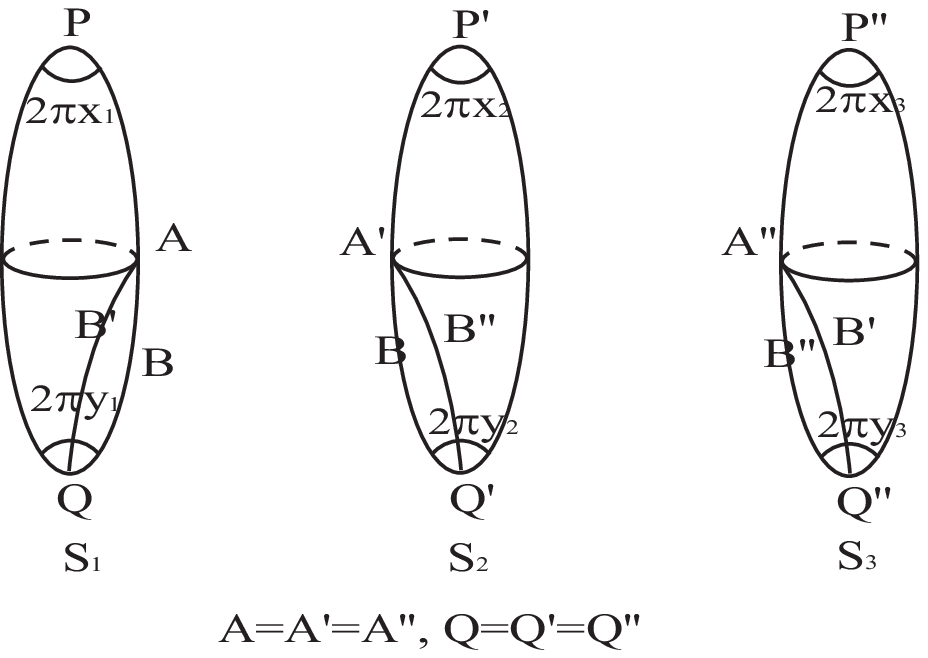}\\
       Figure 3
       \end{center}

     \noindent Along the meridian, cut $S^2_{\{x_1,y_1\}}$ from $A$ to $Q$.
     $\widehat{AQ}$ becomes two identical arcs: $B$ and $B'$. The
     same as above, we cut $S^2_{\{x_2,y_2\}}$ along the meridian
     and get $B,B''$; we cut $S^2_{\{x_3,y_3\}}$ along the
     meridian and get $B',B''$. Then we glue $B$ in
     $S^2_{\{x_1,y_1\}}$ and $B$ in $S^2_{\{x_2,y_2\}}$ together,
     $B'$ in $S^2_{\{x_1,y_1\}}$ and $B'$ in $S^2_{\{x_3,y_3\}}$
     together, $B''$ in $S^2_{\{x_2,y_2\}}$ and $B''$ in
     $S^2_{\{x_1,y_1\}}$ together, to obtain a HCMU metric with
     one saddle point $A\ (=A'=A'')$ of angle $6 \pi$. Meanwhile,
     $Q=Q'=Q''$, at which the angle is $2 \pi(y_1+y_2+y_3)$.\par

Obviously we can glue $\alpha_i$ footballs together in the same
way, to obtain a HCMU metric on $S^2$, which has $\alpha_i$ local
maximum points of angles $2\pi x_1,2\pi x_2,\cdots,2\pi
x_{a_{i}}$, one saddle point of angle $2\pi \alpha_i$, and one
minimum point of angle $2\pi(y_1+y_2+\cdots+y_{\alpha_i})$.

\vskip0.1cm
     \noindent{\it Step 2: construction of a HCMU metric with $j_0$ saddle
     points.} \par
      We choose $\alpha_1$ footballs       $S^2_{\{x_k,y_k\}}~(k=1,2,\cdots,\alpha_1)$
 to construct the first saddle point of angle $2 \pi
      \alpha_1$ like Step 1. Then we choose another meridian on $S^2_{\{x_{\alpha_1},y_{\alpha_1}\}}$ which is different
      from the one passing the previous saddle point, cut this meridian like Step 1, choose  $\alpha_2-1$ footballs
      $S^2_{\{x_k,y_k\}}~(k=\alpha_1+1,\cdots,\alpha_1+\alpha_2-1)$,
      cut them and glue them together with the football $S^2_{\{x_{\alpha_1},y_{\alpha_1}\}}$ like Step 1, we
      get the
      second saddle point of angle $2 \pi \alpha_2$. And then we
      choose $\alpha_3-1$ footballs $S^2_{\{x_k,y_k\}}~(k=\alpha_1+\alpha_2,\cdots,
      \alpha_1+\alpha_2+\alpha_3-2)$
       with
       $S^2_{\{x_{\alpha_1+\alpha_2-1},y_{\alpha_1+\alpha_2-1}\}}$
       to construct the third saddle point of angle $ 2\pi
       \alpha_3$, and so on,  finally we have chosen
$\dsum_{i=1}^{j_0}\alpha_i-(j_0-1)$ footballs to construct a HCMU
metric on $S^2$ with $j_0$ saddle points, as we desire.
\par\vskip0.1cm
       \noindent {\it Step 3: construction of
       $S^2_{\{\alpha_1,\alpha_2,\cdots,\alpha_n\}}$.}\par\vskip0.1cm

       Denote $N=\dsum_{i=1}^{j_0} \alpha_i-(j_0-1)$.
       By Step 2, we have constructed a HCMU metric of $j_0$ saddle
       points with the angles $2 \pi \alpha_1,2 \pi \alpha_2, \cdots,
       2 \pi \alpha_{j_0}$, $N$ maximum points with the angles
        $2 \pi x_k \ (k=1,2,\cdots,N)$, and one minimum point with the angle $2
       \pi \dsum_{k=1}^N y_k$. These angles satisfy:
       \begin{equation}\label{23}\frac{x_k}{y_k}=\frac{x_l}{y_l}>1,~ k,l \in
       \{1,2,\cdots,N\}.\end{equation}
        In the following  we adjust $x_k,\ y_k$ to make the metric
       coincide with what we desire.\par
       Without loss of generality, we assume
       $\displaystyle{\alpha_n=\min_{j_0+1 \le k \le
       n}\{\alpha_k\}}$
        and denote $s=\dsum_{i=1}^{j_0}\alpha_i-(j_0-1)-(n-j_0-1)$.
       By the condition of Theorem A,
       $$s=\sum_{i=1}^{j_0}\alpha_i-n+\chi(S^2)\geq
        0.$$

         {\it Case 1.} If $s=0$, we let
$$ \begin{array}{ll}
x_k=\alpha_{j_0+k} \ \ (1 \le k \le N),\\
\alpha_n=\dsum_{k=1}^N y_k.\end{array}$$
 By (23) the
angles $y_1,y_2\cdots,y_N$ must satisfy the following equations:
\[ \left \{ \begin{array}{l}  \dsum_{k=1}^N y_k =\alpha_n,\\
\dfrac{\alpha_{j_0+1}}{y_1}=\dfrac{\alpha_{j_0+2}}{y_2}=\cdots=\dfrac{\alpha_{n-1}
}{y_{n-j_0-1}}.
\end{array} \right. \]
The equations have a unique solution
$$y_k=\alpha_{j_0+k}{\displaystyle \frac{\alpha_n}{\sum
\limits_{i=j_0+1}^{n-1}\alpha_i}}\ \ ~~(1 \le k \le N=n-j_0-1).
$$
Thus we have proved the claim in this case.\par\vskip0.1cm
  {\it Case 2.}
        If $s>0,$ we let
       \begin{eqnarray*} x_k &=& \alpha_{j_0+k},\ \ (1 \le k \le n-j_0-1),\\
        x_k &=& 1\ \ \ \ \ \ (n-j_0 \le k \le N),\\
        \sum \limits_{k=1}^{N}y_k&=&\alpha_n.
\end{eqnarray*}
This means that there are $s~(=N-(j_0-1))$ smooth maximum points
in the surface, correspondent to angles $x_k ~(k=n-j_0,\cdots, N)$
.

The undetermined angles $y_k~(k=1,\cdots,N)$ satisfy

        \[ \left \{ \begin{array}{l}  \sum \limits_{k=1}^{N}y_k=\alpha_n,\\
\dfrac{\alpha_{j_0+1}}{y_1}=\dfrac{\alpha_{j_0+2}}{y_2}=\cdots=\dfrac{\alpha_{n-1}
}{y_{n-j_0-1}}=\dfrac{1}{y_{n-j_0}}=\cdots=\dfrac{1}{y_{N}}.
\end{array} \right. \]
The equations have a unique solution
\begin{eqnarray*} y_k &=& \alpha_{j_0+k}{\displaystyle \frac{\alpha_n}{s+ \sum
\limits_{i=j_0+1}^{n-1}\alpha_i}}\ \ (1 \le k \le n-j_0-1),\\y_k
&=& {\displaystyle \frac{\alpha_n}{s +\sum
\limits_{i=j_0+1}^{n-1}\alpha_i}}\ \ (n-j_0 \le k \le
N).\end{eqnarray*} We also prove the claim in this case.

Till now we have constructed a HCMU metric on $S^2$ with the
singular angles $\alpha_1, \alpha_2, \cdots, \alpha_n$ and the
correspondent  singular points $q_1,q_2,\cdots,q_n$, which may be
different from the prescribed points $p_1,p_2,\cdots,p_n$.
However, there exists a diffeomorphism $h: S^2 \longrightarrow
S^2_{\{q_1,q_2,\cdots,q_n\}}$, s.t. $h(p_i)=q_i,\ i=1,2,\cdots,n$.
Therefore, we can use $h$ to pull back the HCMU metric $g$
constructed by the claim on $S^2$ to obtain the HCMU metric $h^*g$
as Theorem A desires.$\Box$

\par\vskip0.2cm

 {\bf Remark 2:} We must point out that the HCMU metric satisfying the condition
 in Theorem A is not unique, the metric
 we construct in the proof has exactly one minimum point of the
 scalar curvature $K$. Indeed, by Theorem B we know that the
 construction of a HCMU metric form footballs is a combination problem. If we can
 find some suitable footballs and glue them together to satisfy the condition of
 Theorem A, we have a HCMU metric as desired, and the resulting metric
 may have more than one minimum point of the scalar curvature. For example, we can
 choose two footballs, either $S^2_{\{\frac12,\frac19\}}$ and $S^2_{\{1,\frac29\}}$, or
 $S^2_{\{\frac35,\frac12\}}$ and $S^2_{\{\frac25,\frac13\}}$, to construct
 $S^2_{\{2,\frac12,\frac13\}} $, see following figures for
 illustration.
 \begin{center}

     \includegraphics{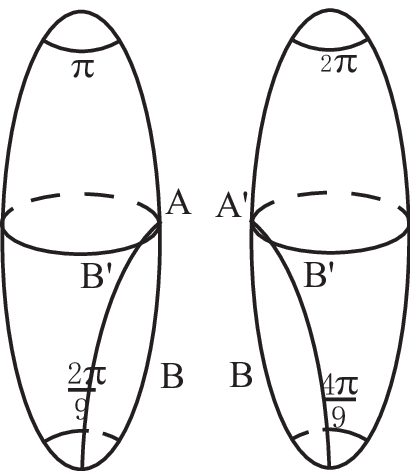} \qquad \qquad
     \includegraphics{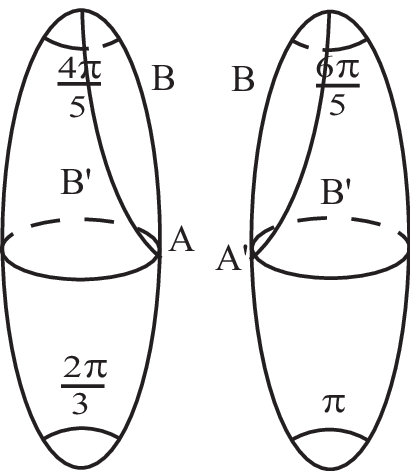}\\
     Figure 4 \qquad \qquad \qquad \qquad \qquad Figure 5
     \end{center}


In fact, Figure 4 is the same as what we constructed in the proof
of Theorem A. However, Figure 5 has two minimum point of the
scalar curvature.

  \par\vskip0.2cm
 {\bf Remark 3:} We should point out that if $\chi(M) \le 0$, the condition (7) is
 not sufficient. For example, if $\chi(M)=-2$, we have the following counter
     example.\par
{\bf Counter example:} Let $\chi(M)=-2,\ ~n=3, ~\alpha_i=2,\ i=1,\
2,\ 3\ \mbox{and} \ j_0=3.$ Then $\sum_{i=1}^{j_0}
\alpha_{\sigma(i)}+\chi(M)-n=6-2-3=1>0,$ but there is no HCMU
metric whose scalar curvature is not a constant s.t. three angles
of saddle points of $K$ are all $4\pi$. If the metric exists,
according to {Proposition 1}, $K$ is continuous in a K-Surface. It
must take its maximum and
  minimum, but the point at which $K$ takes maximum or minimum
  can not be a saddle point, so the point is a smooth critical
  point of $K$. Hence the number of smooth critical
  points of $K$ is more than 1. However,  from the formula $$\chi
  (M)=\sum_{i=1}^j(1-\alpha_i)+n-j+s,\ $$we see $s=1$. That means
  the number of smooth critical points of $K$ is 1, a
  contradiction. $\Box$

  \par\vskip0.2cm
 {\bf Remark 4:} At last we list some problems that might be
 interesting for future study:\par
\begin{enumerate}
\item For other compact Riemannian surfaces, what is the
sufficient and necessary condition for the existence of a HCMU
metric?
\par
\item Is the extremal Hermitian metric unique when none of the
 prescribed angles is an integer multiple of $2\pi$?
\par
\item Given any surface configuration, is the Calabi energy the
only factor determining the connected components
 in the moduli space of HCMU metrics?
\item If we deform the complex structure as well, what is the
structure of the moduli space of HCMU metrics? It will be
interesting to compare this to the classical Teichm\"{u}ller space
in Riemann surfaces.
\end{enumerate}


\par\vskip1cm

\noindent Qing Chen and Yingyi Wu
\\ Department of mathematics\\
University of Science and Technology of China\\
Hefei, Anhui, 230026\\
P. R. China\\
qchen@ustc.edu.cn
\\yingyiw@mail.ustc.edu.cn \\
\vskip0.5cm

\noindent Xiuxiong Chen\\
 Department of Mathematics\\
 University of Wisconsin-Madison\\
 Madison WI 53706\\
 USA\\
 xxchen@math.wisc.edu

\end{document}